\documentclass[12pt]{amsart}
\usepackage{amsmath}

\theoremstyle{plain}
\newtheorem{thm}{Theorem}[section]
\newtheorem{lem}[thm]{Lemma}
\newtheorem{pro}[thm]{Proposition}
\newtheorem{cor}[thm]{Corollary}
\theoremstyle{definition}
\newtheorem{defn}[thm]{Definition}

\newtheorem{rem}[thm]{Remark}

\numberwithin{equation}{section}

\newcommand{\N}{\mathbb{N}}
\newcommand{\R}{\mathbb{R}}

\begin{document}
\title[On a conjecture of Manickam and Singhi]{New results related to a conjecture of Manickam and Singhi}
\date{}
\thanks{This work was partially supported with ``Fondi PRIN (ex40\%)''}



\author{G. Chiaselotti}
\address{Giampiero Chiaselotti, Dipartimento di Matematica, Universit\`{a} della
Calabria, 87036 Arcavacata di Rende, Cosenza, Italy}%
\email{chiaselotti@unical.it}%

\author{G. Infante}
\address{Gennaro Infante, Dipartimento di Matematica, Universit\`{a} della
Calabria, 87036 Arcavacata di Rende, Cosenza, Italy}%
\email{g.infante@unical.it}%

\author{G. Marino}
\address{Giuseppe Marino, Dipartimento di Matematica, Universit\`{a} della
Calabria, 87036 Arcavacata di Rende, Cosenza, Italy}%
\email{gmarino@unical.it}%

\begin{abstract}
In 1998 Manickam and Singhi conjectured that for every positive integer
$d$ and every $n \ge 4d$, every set of $n$ real numbers whose sum is
nonnegative contains at least $\binom {n-1}{d-1}$ subsets of size $d$ whose sums are nonnegative.
In this paper we establish new results related to this conjecture.
We also prove that the conjecture of Manickam and Singhi does not hold for $n=2d+2$.
\end{abstract}
\maketitle
\section{Introduction}
In this paper we establish new results related to a conjecture
of Manickam and Singhi (from now on, (MS)-conjecture).
In order to illustrate the (MS)-conjecture and our results we need to
introduce the following notation.
Let $n\in \N$ and let $I_n$ be the set $\{1,2,\dots,n\}$.
A function $f:I_n\to \R$ is called a $n$--\emph{weight function}
if
$$\sum_{x\in I_n}f(x)\ge0.$$
Let
$W_n(\R)$ denote the set of all $n$--weight functions. If
$f\in W_{n}(\R)$ we set
$$
f^+:=\left|\{x\in I_n:f(x)\geq 0 \}\right|.
$$

If $d$ is an integer with $1\leq d\leq n$ and $Y$ is a subset
of $I_n$ having $d$ elements such that
$$
\sum_{y\in Y}f(y)\ge0,
$$
we call
$Y$ a $(d^+,n)$--\emph{subset} of $f$. If $f \in W_{n}(\R)$, we denote by
$\phi (f,d)$ the number of distinct $(d^{+},n)$--subsets of $f$.

Furthermore, we set
$$
\psi (n,d):=\min\{\phi (f,d): f \in W_{n}(\R)\}.
$$
In 1988, Manickam and Singhi \cite{man-singhi} stated the following conjecture:

\textbf{(MS)-Conjecture }:
If $d$ is a positive integer and $f$ is a $n$--weight function
with $n\ge4d$, then
$$
\psi(n,d)\geq \binom{n-1}{d-1}.
$$
We remark that, as previously observed in \cite{chias}, the
conjecture is equivalent to require that $ \psi(n,d)=
\binom{n-1}{d-1}. $

This conjecture is interesting for several reasons.
It is deeply related with the
first distribution
invariant of the Johnson-scheme \cite{bier-man,man1,man2,man-singhi}. The distribution invariants
were introduced by Bier \cite{bier} and later investigated in \cite{bier-delsarte,man-ths,man1,man-singhi}.
Manickam and Singhi \cite{man-singhi} claim that this conjecture is, in some sense,
dual to the theorem of Erd\"{o}s--Ko--Rado \cite{erd-ko-rad}.
Also, as pointed out by Srinivasan \cite{srini}, this conjecture
settles some cases of another conjecture on multiplicative functions by Alladi, Erd\"{o}s and Vaaler \cite{erdos}.

In general the conjecture of Manickam and Singhi still remains open. So far the following
 partial results have been achieved:
\begin{enumerate}
  \item The (MS)-conjecture is true if $n=ud$, $u\ge4$ (Corollary 1 of \cite{bier-man}).
  \item The (MS)-conjecture is true if $d=2$ (Corollary 3 of \cite{bier-man}).
  \item If we set $n$ in the form $n=ud+v$, where $u\ge4$ and $v=1,\ldots,
d-1$, the (MS)-conjecture is true if $r\le \min \left\{
{\frac{n-v-1}{(d-1)(d-2)}},{\frac{n-v}{d}} \right\}$ (Lemma 2 of \cite{bier-man}).
  \item The (MS)-conjecture is true if $d=3$ and $n\ge93$ (Theorem of \cite{bier-man}).
  \item The (MS)-conjecture is true if $r\le d\le {n/2}$
(Proposition 2 of \cite{chias}).
  \item The (MS)-conjecture is true if it is true when $d<r\le{\frac{d-1}{d}}n$
(Proposition 5 of \cite{chias}).
  \item The (MS)-conjecture is true if $d=3$ (Section 3 of \cite{marino-chias}).
  \item The (MS)-conjecture is true if $n\geq 2^dd^{d+1}+2d^3$ (Theorem 3 of \cite{man-mik}).
  \item The (MS)-conjecture is true if $d>3$ and $n\ge d(d-1)^d(d-2)^d+
d(d-1)^2(d-2)+d[n]_k$, where $[n]_k$ denotes the smallest positive integer congruent to $n$ (mod $k$)
(Main Theorem in \cite{bier-man}).
  \item The (MS)-conjecture is true if $n\geq 2^{d+1}e^{d}d^{d+1}$ (Theorem 1 of \cite{bhatt}).
\end{enumerate}

We point out that for $d>4$ the best estimate between (8), (9) and (10) is (8).

Different techniques have been used to attempt to tackle the (MS)-conjecture.
In \cite{bier-man,chias,man-mik} the approach is combinatorial. In particular Bier and Manickam
\cite{bier-man} use a result of Baranyai (see for example \cite{baranai,wilson}). Manickam and Miklos
\cite{man-mik} use a circle permutation method, previously utilized
by Katona \cite{katona}
for a simpler proof of the theorem of Erd\"{o}s--Ko--Rado.
The approach in \cite{bhatt,bhatt-thesis,marino-chias} is somewhat different. In fact the
techniques in \cite{marino-chias} are analytical-combinatorial and in \cite{bhatt,bhatt-thesis}
are probabilistic.

A natural question arises when one studies the (MS)-conjecture:
\begin{center}
\emph{What is the value of $\psi(n,d)$ for each $d \leq n$?}
\end{center}
In order to provide an answer to this question, in \cite{chias}
the following numbers were introduced:
$$
\gamma(n,d,r)=\min\{\phi(f,d):f\in W_n(\R), f^+=r\},
$$
where $r,d\in \N$, with $r,d\leq n$.

It is clear that a complete computation of these numbers would
also provide a complete determination of the numbers $\psi(n,d)$,
since

\begin{equation}\label{psi-uguale}
  \psi(n,d)=\min_{1\leq r \leq n}\gamma(n,d,r).
\end{equation}

In particular, the knowledge of $\gamma(n,d,r)$ when $n \geq 4d$
and $r$ is an arbitrary integer such that $r \leq n$, would supply
an answer to the (MS)-conjecture.
\begin{rem}\label{chiasgap}
In general the computation of $\gamma(n,d,r)$, started in
\cite{chias}, is not an easy task. For some values of $n,d,r$,
this has been done in \cite{chias}. Nevertheless we stress that
there is a gap in the proof of Proposition 2 of \cite{chias}. In
particular, this means that it is not clear whether the identity
$$
\gamma(n,d,r)=\binom {n-1}{d-1}
$$
holds or not for $r\leq d \leq n/2$.
\end{rem}

In this paper we continue the study of the numbers
$\gamma(n,d,r)$. Here we establish some lower and upper bounds for
$\gamma(n,d,r)$ when $d \leq r \leq \frac{d-1}{d}n$. From these
inequalities we obtain that
\begin{equation}\label{numero-magico}
  \gamma(n,d,r)=\binom {r}{d}+\binom {r}{d-1},
\end{equation}
when $n=2d+2$ with $r=2d-1$ and when
$r=\frac{d-1}{d}n$. Combining our results with the ones in \cite{chias}
we obtain the following values of $\gamma(n,d,r)$:
$$
\gamma(n,d,r)=
\begin{cases}
       \binom{n-1}{d-1}   &{\rm if\ } r\leq d\leq {\frac{n}{2}} \quad ( \star )\cr
       \binom{n-r}{d-r}   &{\rm if\ } r\leq d< n\quad {\rm and\ } \quad r < \frac{n}{n-d}\cr
       \binom{r}{d}       &{\rm if\ } d<r<n     \quad {\rm and\ } \quad r>{\frac{d-1}{d}}n\cr
       \binom{n-1}{d-1}   &{\rm if\ } r=1\cr
       \binom{r}{d} + \binom{r-1}{d-1}  &{\rm if\ } n=2d+2\quad {\rm and\ } \quad r=2d-1\cr
       \binom{r}{d} + \binom{r-1}{d-1}  &{\rm if\ } r\geq d\quad {\rm and\ } \quad r=\frac{d-1}{d}n, \cr
\end{cases}
$$
where $(\star)$ in the first row means that the equality in that
case is uncertain (see Remark \ref{chiasgap}).
We stress that the
determination of the numbers $\gamma(n,d,r)$ in general is an open problem.

A straightforward consequence of \eqref{numero-magico} is that the (MS)-conjecture does not hold if $n=2d+2$.
This provides another range of values of $n$, when $n < 4d$, for which the (MS)-conjecture does not hold.
Note that Bier and Manickam \cite{bier-man} already proved that (MS)-conjecture does
not hold in general. In particular they proved that the (MS)-conjecture does not hold for $n=2d+1$, with $d\geq
2$, and for $n=3d+1$, with $d\geq 3$.

We also prove that, for $n=2d+2$ and $r=2d-1$, \eqref{numero-magico} improves the results of Lemma 1 of
\cite{bier-man}.

A key tool in our paper is Hall's Theorem, as far as we know, used here for the first time in this context. We
use this Theorem to determine, in a non constructive way, certain biunivocal functions between complementary
$q$--subsets of a set with $2q+1$ elements. Such functions are important to compute the numbers $\gamma(n,d,r)$
in the case $n=2d+2$ and $r=2d-1$.

We also suggest a new algorithm to determine the previous functions also in a constructive way.

\section{Preliminaries}\label{preliminari}
In this Section we introduce some notation and prove some elementary arithmetical
preliminaries useful in the sequel of this paper.

We shall assume that a generic weight function
$f \in W_{n}(\R)$, with $f^{+}=r$,
has the form
\begin{equation}\label{forma-f}
   \begin{array}{cccccc}
    1 & \cdots & r & r+1 & \cdots & n \\
    \downarrow & \cdots & \downarrow & \downarrow & \cdots & \downarrow \\
    x_{1} & \cdots &  x_{r} & y_{1} & \cdots & y_{n-r} \\
  \end{array},
\end{equation}
with
$$
x_{1}\geq x_{2}\geq \ldots \geq x_{r}\geq 0 >y_{1}\geq y_{2}\geq \ldots \geq y_{n-r}.
$$
Let us call the indexes $1,\ldots,r$ the \emph{non-negative
elements} of $f$ and the indexes $r+1,\ldots,n$ the \emph{negative
elements} of $f$. The real numbers $x_{1}, \ldots x_{r}$ are said
to be the \emph{non-negative values} of $f$ and the numbers
$y_{1}, \ldots y_{n-r}$ are said to be the \emph{negative values}
of $f$.

If $i_1,\ldots,i_{\alpha}$ are non-negative elements of $f$ and
$j_1,\ldots,j_{\beta}$ are negative elements of $f$, with
$i_1<\ldots<i_{\alpha}$ and $j_1<\ldots<j_{\beta}$, a subset $A$
of $\{1,\ldots,n\}$ is said to be of type
\begin{equation}
  [i_1,\ldots,i_{\alpha}]_{a}^{+}\;[j_1,\ldots,j_{\beta}]_{b}^{-}
\end{equation}
if $A$ is made of $a$ elements chosen in
$\{i_1,\ldots,i_{\alpha}\}$ and $b$ elements chosen in
$\{j_1\ldots,j_{\beta}\}$.

Let $X$ be a finite set of integers. If $q$ is an integer less or
equal than $|X|$, we call $q$--\emph{string} on $X$ a sequence
$a_1\ldots a_q$, where $a_1,\ldots,a_q$ are distinct elements of
$X$ such that $a_1<\ldots <a_q$. The family of all the
$q$--strings on X will be denoted by $X^{(q)}$. In this paper,
each subset $Y$ of $X$ with $q$ elements will be identified with
the $q$--string of his elements ordered in an increasing way. When
$i_1,\ldots,i_k$ are non-negative elements of $f$ and
$j_1,\ldots,j_l$ are negative elements of $f$, with
$i_{1}<\ldots<i_{k}<j_{1}<\ldots<j_{l}$, the $(k+l)$--string
$i_1\ldots i_k j_1\ldots j_l$ will be written in the form
$$
i_{1}\ldots i_{k}|(j_{1}-r)\ldots(j_{l}-r),
$$
 (thus $j_{1}-r,\ldots,j_{l}-r\in \{1,\ldots,n-r\}$).
\newline

For example, if $n=10$ and $r=7$, the 4--string $1269$ will be
written in the form $126|2$.
\newline

Using the string-terminology instead of the set-terminology, in the sequel we call a $(d^+,n)$--subset of $f$
a $(d^+,n)$--\emph{string} of $f$.
\newline

Let us consider now the partition $\mathcal{P}$ of the real
interval $(0,\frac{d-1}{d}n]$:
\begin{equation}\label{partizione}
\mathcal{P}=\Bigl\{0, \frac{d-1}{d}, 2\frac{d-1}{d}, \ldots,(n-1)\frac{d-1}{d},n\frac{d-1}{d} \Bigr\}.
\end{equation}
The following Proposition establishes when an interval determined by $\mathcal{P}$ contains an integer.
\begin{pro}\label{propA}
If $k=1,\ldots,n$ and if $n-k\not\equiv_{d} 0$, there exists a unique integer $r$ such that
\begin{equation}\label{range-r2}
\frac{d-1}{d}(n-k)<r\leq \frac{d-1}{d}(n-k+1),
\end{equation}
and $r$ coincides with $\lfloor \frac{d-1}{d}(n-k+1) \rfloor$. Furthermore if $n-k\equiv_{d} 0$
no integer $r$ satisfies \eqref{range-r2}.
\end{pro}
\begin{proof}
Let $k \in \{1,\ldots,n\}$ and set $m=n-k+1$. Since the interval $(\frac{d-1}{d}(m-1), \frac{d-1}{d}m]$ has
length $\frac{d-1}{d}<1$, there is at most one integer $r$ that satisfies \eqref{range-r2}.
Let us now write $m$ in the form
\begin{equation}\label{forma-m}
m=\tilde{q}d+s,
\end{equation}
where $\tilde{q},s$ are integers such that
$\tilde{q}\geq 0$, $1\leq s \leq d$.
Let us suppose now that $n-k \not \equiv_{d} 0$, that is $m\not \equiv_{d} 1$; then we have $2\leq s \leq d$.\newline
Let $r=\lfloor \frac{d-1}{d}m \rfloor$. We show that $r$ satisfies \eqref{range-r2}.\newline
Firstly, the second inequality of \eqref{range-r2} is straightforward; secondly, for the first inequality we observe
$$
\frac{d-1}{d}(m-1)=\frac{d-1}{d}(\tilde{q}d+s-1)=\tilde{q}(d-1)+(s-1)\frac{d-1}{d}.
$$
Furthermore
\begin{equation}\label{forma-m1}
r=\lfloor \frac{d-1}{d}(\tilde{q}d+s) \rfloor=\lfloor \tilde{q}(d-1)+s -\frac{s}{d}\rfloor
=\tilde{q}(d-1)+(s-1).
\end{equation}
Therefore $r>\frac{d-1}{d}(m-1)$, since $s\geq 2$.

If $n-k\equiv_{d} 0$, that is $m \equiv_{d} 1$, in \eqref{forma-m} we have $s=1$ and \eqref{range-r2} becomes
\begin{equation}\label{stella}
\tilde{q}(d-1)<r\leq \tilde{q}(d-1)+\frac{d-1}{d}.
\end{equation}
Note that \eqref{stella} has no integer solutions.
\end{proof}
\begin{lem}\label{propB}
Let $r$ be a positive integer such that
\begin{equation}\label{range-r}
d\leq r\leq \frac{d-1}{d}n.
\end{equation}
Then there exists a unique positive integer $b(r) \in
\{1,\ldots,n-r-1\}$ that satisfies
\begin{equation}\label{slup}
\frac{d-1}{d}(n-b(r))<r\leq \frac{d-1}{d}(n-b(r)+1),
\end{equation}
\end{lem}
\begin{proof}
By construction of partition $\mathcal{P}$, as in
\eqref{partizione}, there exists a unique $b(r) \in
\{1,\ldots,n\}$ such that \eqref{slup} holds.

We now show that $b(r)$ cannot exceed $n-r-1$.

Firstly, we suppose that $b(r)>n-r$. Then, we write $b(r)$ in the
form $b(r)=n-r+\zeta$, with $\zeta$ integer such that $1\leq \zeta
\leq r$. Since $r$ satisfies \eqref{slup}, we have
$$
\frac{d-1}{d}(r-\zeta)<r\leq \frac{d-1}{d}(r-\zeta+1),
$$
that is
\begin{equation}\label{range-r3}
    \zeta (1-d)<r\leq (d-1) (1-\zeta).
\end{equation}
Since $\zeta (1-d)<0$ and $(d-1) (1-\zeta)\leq 0$, there is no positive integer $r$
that satisfies \eqref{range-r3}.

Secondly, if $b(r)=n-r$, \eqref{slup} becomes
$$
\frac{d-1}{d}r<r\leq \frac{d-1}{d}(r+1),
$$
that is
$$
0<r\leq d-1,
$$
contradicting the hypothesis \eqref{range-r}.
\end{proof}
We stress that the number $b(r)$ will play a key role in the
sequel of the paper.
\section{Some upper and lower bounds for $\gamma(n,d,r)$}\label{stime-gamma}
In this Section we establish an upper bound for $\gamma(n,d,r)$, when $r$ satisfies \eqref{range-r}.
We also provide a lower bound for $\gamma(n,d,r)$ under one additional hypothesis.
\begin{pro}\label{propC}
Let $r$ be a positive integer that satisfies 
$$
d\leq r \leq \frac{d-1}{d}n,
$$
then
\begin{equation}\label{gamma-UB}
\gamma(n,d,r)\leq \sum_{j=0}^{\min\{b(r),\, d-1\}} \binom
{b(r)}{j}\binom {r}{d-j}.
\end{equation}{}
\end{pro}
\begin{proof}
Since $1\leq b(r) \leq n-r-1$, we construct a weight function
$f\in W_{n}(\R)$, with $f^{+}=r$, such that
\begin{equation}\label{somma-phi}
\phi(f,d)= \sum_{j=0}^{\min\{b(r),\, d-1\}} \binom {b(r)}{j}\binom
{r}{d-j}.
\end{equation}
This is sufficient to prove the thesis.

Let $h=\min\{b(r),\, d-1\}$. Let $\alpha$ be a positive real
number. In order to simplify the notation, we call $\beta$ the number
$\frac{r+b(r)(-\alpha)}{n-r-b(r)}$, in such a way that
$$
r+b(r)(-\alpha)+(n-r-b(r))(-\beta)=0
$$
holds.

At this point we define the function
\begin{equation}\label{forma-f-bast}
f_{\alpha}:
   \begin{array}{ccccccccc}
    1 & \cdots & r & r+1 & \cdots & r+b(r) & r+(b(r)+1) & \cdots & r+(n-r)\\
    \downarrow & \cdots & \downarrow & \downarrow & \cdots & \downarrow & \downarrow & \cdots & \downarrow\\
    1 & \cdots &  1 & -\alpha & \cdots & -\alpha & -\beta  & \cdots & -\beta\\
  \end{array}.
\end{equation}
We now show that for $\alpha$ sufficiently small, that is
\begin{equation}\label{alpha}
0<\alpha<\min\Bigl\{\frac{r}{b(r)},\; \frac{d}{h}-1,\; \frac{d}{b(r)}\bigl(r-\frac{d-1}{d}(n-b(r))\bigr)\Bigr\},
\end{equation}
$f_{\alpha}$ is a weight function that satisfies \eqref{somma-phi}.

In fact:
\begin{itemize}
    \item[$a)$] the denominator of $\beta$, due to Lemma \ref{propB}, is a positive number. Furthermore the numerator
of $\beta$ is a positive number if and only if $\alpha<\frac{r}{b(r)}$.
Therefore \eqref{alpha} and the definition of $\beta$ assure that $f_{\alpha}$ is a weight function.
    \item[$b)$] having $\alpha < d/h-1$ is equivalent to require
\begin{equation}\label{condB}
\underbrace{1+\ldots+1}_{d-h \text{ times} }+\underbrace{(-\alpha)+\ldots+(-\alpha)}_{h \text{ times} }\geq 0.
\end{equation}
This condition assures that
the subsets of the type
\begin{equation}\label{condsets}
  \begin{array}{lccc}
   [1,\ldots,r]_{d}^{+} & [r+1,\ldots,r+b(r)]_{0}^{-}, & \text{in total} & \binom {b(r)}{0}\binom {r}{d}\\
   {} [1,\ldots,r]_{d-1}^{+} & [r+1,\ldots,r+b(r)]_{1}^{-}, & \text{in total} & \binom {b(r)}{1}\binom {r}{d-1} \\
   \phantom{aaaa}\vdots & \vdots & \vdots & \vdots\\
   {} [1,\ldots,r]_{d-h}^{+} & [r+1,\ldots,r+b(r)]_{h}^{-}, & \text{in total} & \binom {b(r)}{h}\binom {r}{d-h} \\
  \end{array}
\end{equation}
are $(d^+,n)$--subsets of $f$
    \item[$c)$] firstly we note that the requirement
    $$\alpha < \frac{d}{b(r)}\Bigl(r-\frac{d-1}{d}(n-b(r))\Bigr)$$
    is equivalent to require
    \begin{equation}\label{condC}
\frac{d-1}{d}(n-b(r))+\alpha \frac{b(r)}{d}<r.
\end{equation}
Lemma \ref{propB} assures the existence of a such $\alpha$. Note that
\eqref{condC} is equivalent to
   \begin{equation}\label{condA}
\underbrace{1+\ldots+1}_{d-1 \text{ times} }+(-\beta)<0,
\end{equation}{}
that assures that the $(d^+,n)$--strings of $f$ are only of the type \eqref{condsets}.
Therefore we have constructed a weight function $f$ with $r$ non-negative elements that satisfies
\eqref{somma-phi}.
\end{itemize}
\end{proof}

\begin{cor}\label{Cor1}
Let $r$ be a positive integer such that $r\geq d$ and $\frac{d-1}{d}(n-1)< r \leq \frac{d-1}{d}n$. Then
\begin{equation}\label{upper-gamma}
\gamma(n,d,r)\leq \binom {r}{d}+\binom {r}{d-1}.
\end{equation}
\end{cor}
\begin{proof}
The result follows directly from Proposition \ref{propC} since $b(r)=1$.
\end{proof}

\begin{pro}\label{propD}
Let r a positive integer such that $d \leq r \leq \frac{d-1}{d}n$. Let $f \in W_{n}(\R)$,
with $f^{+}=r$, as in \eqref{forma-f}.
 If
\begin{equation}\label{dis1}
x_{1}+y_{n-r}\geq 0,
\end{equation}
then
\begin{equation}\label{dis2}
\phi(f,d)\geq \binom {r-1}{d-2}(n-r)+\binom {r}{d} \geq \binom {r}{d}+\binom {r}{d-1}.
\end{equation}
\end{pro}
\begin{proof}
We can consider
the $d$--strings of $ \{1,\ldots,n\}$ of type
\begin{equation}\label{condD}
1i_{1}\ldots i_{d-2}|(n-r),
\end{equation}
where $i_{1}\ldots i_{d-2}$ are chosen in $\{2,\ldots,r\}$.

By virtue of \eqref{dis1}, each string of the type \eqref{condD} is a
$(d^+,n)$--string of $f$.

On the other hand, since $y_{1}\geq
y_{2}\geq \ldots \geq y_{n-r}$,  each string of type
\begin{equation}\label{condE}
1i_{1}\ldots i_{d-2}|k,
\end{equation}
where $i_{1}\ldots i_{d-2}$ are chosen in $\{2,\ldots,r\}$ and
 $k$ in $\{1,\ldots,n-r\}$, will be a
$(d^+,n)$--string of $f$.

The distinct strings of the type \eqref{condE} are exactly $\binom
{r-1}{d-2}(n-r)$. There are moreover all the
$(d^+,n)$--strings of $f$ that are the $d$--strings on
$\{1,\ldots,r\}$. This proves the first inequality in \eqref{dis2}.
Moreover, since $r\leq \frac{d-1}{d}n$, we also have $n-r \geq \frac{r}{d-1}$.
Therefore
$$
\binom {r-1}{d-2}(n-r)\geq  \binom {r-1}{d-2}\frac{r}{d-1}=\binom {r}{d-1}.
$$
Thus the second inequality also holds.
\end{proof}

As a direct consequence of Corollary \ref{Cor1} and Proposition \ref{propD} it follows that if
$r$ is a positive integer with $r\geq d$ such that $\frac{d-1}{d}(n-1)< r \leq \frac{d-1}{d}n$,
then
\begin{equation}\label{min1}
\min \{ \phi (f,d): f\in W_{n}(\R), f^{+}=r, x_1+y_{n-r} \geq 0\}=\binom {r}{d}+\binom {r}{d-1}.
\end{equation}
\begin{rem}\label{oss}
We conjecture that
\begin{equation}\label{magic-conjecture}
\gamma(n,d,r)= \binom {r}{d}+\binom {r}{d-1},
\end{equation}
when $r\geq d$ and $\frac{d-1}{d}(n-1)< r \leq \frac{d-1}{d}n$.

In Section 5 we give a partial answer to this conjecture. Note that, in order to prove
\eqref{magic-conjecture},
 by Corollary \ref{Cor1}
it is sufficient to show
\begin{equation}\label{lower-gamma}
\gamma(n,d,r)\geq \binom {r}{d}+\binom {r}{d-1}.
\end{equation}
Moreover, by virtue of \eqref{min1}, the inequality \eqref{lower-gamma} is equivalent to the following:

\begin{multline}\label{min2}
\min \{ \phi (f,d): f\in W_{n}(\R), f^{+}=r, x_k+y_{n-r}< 0,\\
 \text{for every}\ k=1,\ldots,r\} \geq \binom {r}{d}+\binom {r}{d-1}.
\end{multline}
In Section 5 we shall prove this inequality in the special case $n=2d+2$.
\end{rem}
We close this section providing a simple combinatorial interpretation of the inequalities

$$\frac{d-1}{d}(n-1)< r \leq \frac{d-1}{d}n. $$

For this purpose let us note that the last inequalities are equivalent to the following:
\begin{equation}\label{rangeComb-r}
(n-r-1)(d-1)<r\leq (n-r)(d-1),
\end{equation}

Let now $r$ be a positive integer that satisfies \eqref{rangeComb-r} and
 $f\in W_{n}(\R)$, with $f^{+}=r$, as in \eqref{forma-f}.
Let us consider the following representation
\begin{equation}\label{rappresentazione}
\begin{array}{rcl}
  \llcorner\lrcorner   + \llcorner\lrcorner + & \ldots &+ \llcorner\lrcorner + k_{1}\\
  \llcorner\lrcorner   + \llcorner\lrcorner + & \ldots &+ \llcorner\lrcorner + k_{2}\\
  &\ldots & \\
  \llcorner\lrcorner   + \llcorner\lrcorner + & \ldots &+ \llcorner\lrcorner + k_{n-r-1}\\
  \llcorner\lrcorner   + \llcorner\lrcorner + &\ldots &+ \llcorner\lrcorner + k_{n-r}\\
  \end{array},
\end{equation}
where every $\llcorner\lrcorner$ can be seen as a ``box'' initially empty and every row contains $d-1$ boxes.
Every of such boxes
can be occupied by at most one non-negative element of $f$. Thus \eqref{rangeComb-r} is equivalent
to state that $n-r-1$ rows in \eqref{rappresentazione} must be completely occupied, whereas
the last row must contain at least a non-empty box and, furthermore, the number of non-negative
elements of $f$ cannot exceed the number of empty boxes in \eqref{rappresentazione}. This combinatorial
interpretation of \eqref{rangeComb-r} suggests to examine firstly the
$(d^{+},n)$--strings of $f$ of the form $+\ldots+-$, that is a subset with $d-1$ non-negative elements
and only one negative.

\section{An application of Hall's Theorem}\label{appl-Hall}
In this Section we use Hall's theorem on distinct
representatives to determine some biunivocal functions between
$q$--subsets of a set with $2q+1$ elements.
The results of this
Section are used in Section \ref{n=2d+2} to determine
$\gamma(n,d,r)$ when $r= 2d-1$ and $n=2d+2$.

We now introduce some definitions and notation useful in the
sequel.

Let $\Omega =
\{1,2,\ldots,2q,2q+1\}$, where $q$ is a fixed positive integer.

Given a $q$--string $a_1\ldots a_q \in \Omega^{(q)}$, for notation convenience we denote by
$\mathfrak{C}_{q}(a_1 \ldots a_q)$
the family of all the $q$--strings
on $ \Omega \setminus \{a_1,\ldots,a_q\}$, that is
$$
\mathfrak{C}_{q}(a_1 \ldots a_q) = (\Omega \setminus \{a_1,\ldots ,a_q\}) ^{(q)}
= \{b_1\ldots b_q : b_1,\ldots ,b_q \in \Omega, b_i \neq a_j, i,j = 1,\ldots ,q\}.
$$
Note that the family $\mathfrak{C}_{q}(a_1 \ldots a_q)$ has exactly $\binom{q+1}{q} = q+1$
distinct $q$--strings.

A $q$--string in $\mathfrak{C}_{q}(a_1 \ldots a_q)$ will
be called a $q$--\emph{almost-complementary} (or $q$-AC) of $a_1,\ldots a_q$.

From now on we call $p$ the number of the distinct $q$--strings of $\Omega^{(q)}$, that is $p=\binom{2q+1}{q}$.
We denote by $A_1,\ldots ,A_p$ all the $q$--strings of $\Omega^{(q)}$ such that
$$
A_1 \prec \ldots \prec A_p,
$$
where $\prec$ is the usual lexicographic order.

\begin{defn}\label{defI}
A $q$--\emph{pairing of almost-complementaries on $\Omega$} (or $q$-PAC on $\Omega$) is a biunivocal function $\varphi :
\Omega^{(q)}\rightarrow \Omega^{(q)}$ such that $\varphi(A_i)$ is a $q$-AC of $A_i$ for $i=1,\ldots p$,
that is
$$
\varphi(A_i) \in \mathfrak{C}_{q}(A_i),
$$
for $i=1,\ldots, p$.
\end{defn}

Let us set now $\mathfrak{F}_{q}=\{\mathfrak{C}_{q}(A_1),\ldots, \mathfrak{C}_{q}(A_p)\}$.

We recall that the family $\mathfrak{F}_{q}$ has a system of distinct representatives (SDR),  say $(C_1,\ldots,C_p)$,
if $C_1\in \mathfrak{C}_{q}(A_1),\ldots, C_p\in \mathfrak{C}_{q}(A_p)$ and $C_i \neq C_j$ for $i,j \in \{1,\ldots,p\}$
with $i\neq j$.
\begin{pro}\label{propE}
The family $\mathfrak{F}_{q}$ has a SDR if and only if there exists a $q$--PAC on $\Omega$.
\end{pro}
\begin{proof}
Sufficiency. Let $(C_1,\ldots, C_p)$ be a SDR for $\mathfrak{F}_{q}$.
This means that all the $C_i$ are distinct $q$--strings and that $C_i \in \mathfrak{C}_{q}(A_i)$ for $i=1, \ldots, p$.
Thus the function
$$
\varphi : \Omega^{(q)} \to \Omega^{(q)}
$$
defined by
$$
\varphi (A_i)=C_i \in \mathfrak{C}_{q}(A_i),
$$
for $i=1,\ldots , p$, is a $q$-PAC on $\Omega$.

Necessity. If $\varphi$  is a $q$--PAC on $\Omega$, then $\varphi$ is a bijection such that
$\varphi (A_i)\in \mathfrak{C}_{q}(A_i)$, for $i=1,\ldots , p$.
Since $\varphi$ is a bijection, $\varphi(A_1),\dots,\varphi(A_p)$ is a SDR for $\mathfrak{F}_{q}$.
\end{proof}
\begin{pro}\label{propF}
For every positive integer $q$ there exists a $q$--PAC on $\Omega$.
\end{pro}
\begin{proof}
By virtue of Proposition \ref{propE} it is sufficient to show that the family $\mathfrak{F}_{q}$ has a
SDR, i.e. that the well-known Hall's condition holds:
\begin{equation}\label{Hall}
\text{for every}\;I\subset \{1,\ldots ,p\},\; \text{we have} \; \Bigl |\bigcup_{i\in I}\mathfrak{C}_{q}(A_{i}) \Bigr |\geq|I|.
\end{equation}
Therefore, let $I=\{i_1,\ldots,i_k\}$ be an arbitrary subset of indices
$i_1,\ldots,i_k \in \{1,\ldots ,p\}$.
Let $Y:=\mathfrak{C}_{q}(A_{i_1})\cup \dots \cup \mathfrak{C}_{q}(A_{i_1k})=\{C_1,\dots,C_a\}$.
With this notation \eqref{Hall} is equivalent to $a \geq k$, therefore we shall prove now this last inequality.
Set $\mathfrak{A}:=\{ \mathfrak{C}_{q}(A_{i_1}), \dots , \mathfrak{C}_{q}(A_{i_k})\}$.
For all $C_l\in Y$ we denote by $d_{\mathfrak{A}}(C_l)$ the degree of $C_l$ respect to the family
$\mathfrak{A}$, that is the number of distinct sets $ \mathfrak{C}_{q}(A_{i_j})$ that contain $C_l$.
We have previously observed that $|\mathfrak{C}_{q}(A_{i_j})|=q+1$ for all $A_{i_j}$, moreover, by a classical double
counting principle we also have
$$
\sum_{l=1}^{a}d_{\mathfrak{A}}(C_l)=\sum_{j=1}^{k}| \mathfrak{C}_{q}(A_{i_j})|,
$$
hence
\begin{equation}\label{eqDEGREE}
\sum_{l=1}^{a}d_{\mathfrak{A}}(C_l)=k(q+1).
\end{equation}

On the other hand, every $C_l$ is a $q$--string, let us say $C_l = c_1\ldots c_q$, which belongs to the sets
$\mathfrak{C}_{q}(a_1\ldots a_q)$, where $a_1\ldots a_q$ is a $q$-AC of $c_1\ldots c_q$.
Since the number of the distinct $q$-AC strings of $c_1\ldots c_q$ is $\binom{q+1}{q}$, it follows that
every $C_l$ belongs exactly to $q+1$ subsets $\mathfrak{C}_{q}(A_{s})$, with $s=1,\ldots,p$;
therefore
$d_{\mathfrak{A}}(C_l)\leq q+1$ for $l=1, \ldots , a$. By \eqref{eqDEGREE}
we obtain then
$$
k(q+1)\leq a (q+1),
$$
i.e. $k\leq a $.
\end{proof}
The Proposition \ref{propF} does not provide an explicit construction of a $q$--PAC on $\Omega$.
In order to construct a $q$--PAC on $\Omega$ we suggest the following Algorithm:
\medskip
\begin{center}
\textbf{$q$--PAC Algorithm}
\end{center}

\noindent  \underline{Input:} a positive integer $q$

\noindent  \underline{Output:} a $q$--PAC on $\Omega$

\begin{enumerate}
  \item[\textit{Step 1:}] Write all the $q$--strings of $\Omega^{(q)}$ ordered in increasing
  way with respect to the lexicographic order
  $$
  B_1\prec B_2\prec\ldots\prec B_p,
  $$
  and put them in an array \verb"Dom[p]" of $q$--strings, that has $p$ positions, where $p=\binom {2q+1}{q}$.
  \item[\textit{Step 2:}] For all $i=1,\ldots,p$ write all the $q$--strings of $\mathfrak{C}_{q}(B_{i})$
  in decreasing lexicographic order
 $$
  C_{i_{1}}\succ\ldots \succ C_{i_{q+1}}.
  $$
  \item[\textit{Step 3:}] Set up an array \verb"Im[p]" of $q$--strings, that has $p$ positions, and initialize
  every position with the $q$--string with all zero entries.
  \item[\textit{Step 4:}] For all $i=1,\ldots,p$ examine in sequence the $q$--strings $C_{i_{1}},\ldots,C_{i_{q+1}}$ and
  put the first of such $q$--strings that does not appear in \verb"Im[1]"$,\ldots,$\verb"Im[i]" in the position \verb"Im[i]".
\end{enumerate}

Then the correspondence \verb"Dom[i]"$\longmapsto$\verb"Im[i]" $(i=1,\ldots, p)$ provides a $q$--PAC on $\Omega$.
\newline

For small values of $q$ we have implemented the previous algorithm in Java. For example, if $q=3$ then
$p=\binom{2q+1}{q} = \binom{7}{3} = 35$. In this case \verb"Dom[35]" and \verb"Im[35]" are two arrays with 35
position, both containing all the 3--strings on $\{1,\ldots,7\}$. The execution of our program for $q=3$
provides the following result (the strings on the left of \verb"--->" are those of \verb"Dom[35]", the strings
on the right of \verb"--->" are those of \verb"Im[35]"):
\begin{quote}
\begin{tiny}
\begin{flushright}
\verb"123 ---> 567; 124 ---> 367; 125 ---> 467; 126 ---> 457; 127 ---> 456;"\\
              \verb"134 ---> 267; 135 ---> 247; 136 ---> 257; 137 ---> 256;"\\
                            \verb"145 ---> 237; 146 ---> 357; 147 ---> 356;"\\
                                          \verb"156 ---> 347; 157 ---> 346;"\\
                                                        \verb"167 ---> 345;"

              \verb"234 ---> 167; 235 ---> 147; 236 ---> 157; 237 ---> 156;"\\
                            \verb"245 ---> 137; 246 ---> 135; 247 ---> 136;"\\
                                          \verb"256 ---> 134; 257 ---> 146;"\\
                                                        \verb"267 ---> 145;"

                            \verb"345 ---> 127; 346 ---> 125; 347 ---> 126;"\\
                                          \verb"356 ---> 124; 357 ---> 246;"\\
                                                        \verb"367 ---> 245;"\\

                                          \verb"456 ---> 123; 457 ---> 236;"\\
                                                        \verb"467 ---> 235;"\\

                                                        \verb"567 ---> 234;"
\end{flushright}
\end{tiny}
\end{quote}
\section{The case $d\leq r$, $\frac{d-1}{d}(n-1)< r \leq \frac{d-1}{d}n$, $n=2d+2$}\label{n=2d+2}

In this Section we shall assume that $n=2d+2$ and that $r$ is a positive integer such that $r\geq d$,
$\frac{d-1}{d}(n-1)< r \leq \frac{d-1}{d}n$. Under such hypotheses we can apply the Proposition \ref{propA}
to the case $k=1$, obtaining
$$r=\lfloor \frac{d-1}{d}n\rfloor=\lfloor \frac{d-1}{d}(2d+2)\rfloor=2d-1=n-3.$$

For such values of $r$ and $n$ we determine the value of $\gamma (n,d,r)$.
This result implies that in this case the (MS)-conjecture does not hold. We also compare our
results with the ones in \cite{bier-man}.

\begin{thm}\label{teorema1}
If $n=2d+2$ and $r=2d-1=n-3$ then
\begin{equation}
\gamma(n,d,r)= \binom {r}{d}+\binom {r}{d-1}.
\end{equation}
\end{thm}
\begin{proof}
Due to Remark \ref{oss} we only need to show \eqref{min2} when $n-r=3$. Thus
take $f\in W_{n}(\R)$, with $f^{+}=r$, as in \eqref{forma-f} and suppose that
$x_k+y_3< 0$ for every $k=1,\ldots,r$.

Take $q=d-1$ (and therefore $r=2q+1$). By Proposition \ref{propF} there exists a $q$--PAC
on $\Omega$,
where $\Omega=\{1,\ldots,r\}=\{1,\ldots,2q+1\}$.
We use the notation introduced in Section \ref{appl-Hall}.
Take $\Omega^{(q)}=\{A_{1},\ldots,A_{p}\}$ with the lexicographic order:
$$
A_1 \prec \ldots \prec A_p,
$$
where $p=\binom{2q+1}{q}$.
Let $C_s=\varphi (A_s)$ for $s=1,\ldots,p$.

Since $A_s$ and $C_s$ are $q$--strings with no common elements there exists in $\Omega$ a unique element,
say $i_s$, that is not an element of the $q$--string $A_s$ and nor an element of the $q$--string $C_s$.
We point out that the elements $i_1,\ldots,i_p$ are not distinct between them, since $p>r$.

If $A= t_1 \ldots t_{d-1} \in \Omega^{(q)}$ and $k \in \{1,2\}$, with the notation $A|k$ we mean the
$d$--string $t_1 \ldots t_{d-1}|k$ and with $i_{s}|3$ the $2$--string
with the non-negative element $i_{s}$ and with the negative element $3$.
We now consider the following configuration:
\begin{equation}\label{configuration}
\begin{array}{rcl}
  A_1|1 & C_1|2 & i_1|3\\
  A_2|1 & C_2|2 & i_2|3\\
  &\ldots & \\
  A_p|1 & C_p|2 & i_p|3\\
  \end{array}.
\end{equation}
Since $\varphi$ is a bijection, the $q$--strings $C_1|2,\ldots,C_p|2$ are themselves
distinct.
Moreover, since $\varphi$ is a $q$--PAC
on $\Omega$, each row in \eqref{configuration} contains all the elements (non-negative and negative) of $f$.
Since the function $f$ is a weight function and from the hypothesis we have $x_{i_{s}}+y_{3}< 0$ (that is
 each $i_s|3$ corresponds to a negative sum), in every $s$th-row of the
configuration \eqref{configuration} at least one $d$--string between
$A_s|1$ and $C_s|2$ must be a $(d^+,n)$--string for $f$.

This shows that the number of the distinct $(d^+,n)$--strings for $f$ is at least
equal to the number of the rows in \eqref{configuration}, that is $p=\binom{2q+1}{q}=\binom{r}{d-1}$.

The remaining $(d^+,n)$--strings for $f$ that we need in order to obtain \eqref{min2} are all the $d$--strings of $\Omega$,
which are $\binom{r}{d}$. This shows that $\phi(f,d) \geq \binom{r}{d} + \binom{r}{d-1}$.

Since $f$ is arbitrary, \eqref{min2} follows.
\end{proof}

From this result we deduce the following consequence on the (MS)-conjecture.
\begin{cor}
  The (MS)--conjecture does not hold when $n=2d+2$ and $d\geq 3$.
\end{cor}
\begin{proof}
We take $r=2d-1$. Then from \eqref{psi-uguale} and Theorem \ref{teorema1} we have
  $$
    \psi(n,d) \leq \binom{r}{d} + \binom{r}{d-1},
  $$
and for such values of $r$ and $n$ we have
$$
\binom{r}{d} + \binom{r}{d-1} = \frac{2(2d-1)!}{d!(d-1)!} < \frac{(2d+1)!}{(d-1)!(d+2)!} = \binom{n-1}{d-1},
$$
if $d\geq 3$.
\end{proof}

When $n=2d+2$ and $r=2d-1$,
from Theorem \ref{teorema1} it follows that
\begin{align*}
\phi(f,d) \geq \binom{2d-1}{d} + \binom{2d-1}{d-1} &>
\binom{2d-1}{d-1}= \binom{(2d+2)-2-1}{d-1},
\end{align*}
that is
\begin{equation}
\phi(f,d)>\binom{n-k-1}{d-1}.
\end{equation}
This inequality improves the estimate
\begin{equation}\label{BM}
    \phi(f,d)\geq \binom{n-k-1}{d-1},
\end{equation}
obtained in \cite{bier-man} under the additional hypotheses
\begin{enumerate}
  \item[$(i)$] $\sum_{x \in I_n} f(x) = 0$,
  \item[$(ii)$] $\sum_{y \in Y} f(y) \neq 0$ for all $Y \subseteq I_n$ such that $|Y|=d$.
\end{enumerate}
\section{The case $r=\frac{d-1}{d}n$}\label{rd=(d-1)n}
There is also another case when we can prove
$$
\gamma(n,d,r)=\binom {r}{d}+\binom {r}{d-1}.
$$
This is the case $r=\frac{d-1}{d}n$. In such case $d|n$ and therefore $r\geq d$ if $d\geq 2$.
\begin{thm}\label{teorema2}
Let $r=\frac{d-1}{d}n$. Then
$$
\gamma(n,d,r)=\binom {r}{d}+\binom {r}{d-1}.
$$
\end{thm}
\begin{proof}
The condition $r=\frac{d-1}{d}n$ is equivalent to $r=(d-1)(n-r)$.
Take $f\in W_n(\R)$ with $f^+ = r$. Then we can build partitions
$\mathfrak{S}$ of the set $\{1,\ldots,r,r+1,\ldots,n\}$ of the
type $\mathfrak{S}=\{C_1,\ldots,C_{n-r}\}$, where
\begin{equation}\label{BM-part}
\begin{array}{rl}
  C_1 & =i_{1,1} \ldots i_{1,d-1}|k_{1}\\
  C_2 & =i_{2,1} \ldots i_{2,d-1}|k_{2}\\
  &\ldots \\
  C_{n-r} & =i_{n-r,1} \ldots i_{n-r,d-1}|k_{n-r}\\
  \end{array},
\end{equation}
with $i_{s,t} \in \{1,\ldots , r\}$, $k_{l}\in \{1,\ldots, n-r\}$,
for $1\leq s \leq n-r$, $1\leq t \leq d-1$, $1\leq l \leq n-r$. By
means of a technique similar to the one used by Bier and Manickam
in the proof of Lemma 1 of \cite{bier-man}, we can claim that
there exist exactly $\binom {r-1}{d-2}(n-r)$ disjoint partition of
type \eqref{BM-part}. Since $f$ is a weight function, at least a
row in \eqref{BM-part} is $(d^+,n)$--string for $f$.

Since the partitions are disjoint, if we extract from each of them
at least a $(d^+,n)$--string for $f$, we get at least  $\binom
{r-1}{d-2}(n-r)$ distinct $(d^+,n)$--string for $f$.

Since $r=\frac{d-1}{d}n$, we have
$$
\binom {r-1}{d-2}(n-r)=\binom {r-1}{d-2}\frac{r}{d-1}=\binom
{r}{d-1}.
$$
Therefore
$$
\phi(f,d) \geq \binom{r}{d} + \binom{r}{d-1}.
$$
This show that
$$
\gamma(n,d,r) \geq \binom{r}{d} + \binom{r}{d-1}.
$$
The equality follows form Corollary \eqref{Cor1}.
\end{proof}

\end{document}